\documentclass[12pt]{article}
\usepackage{amsfonts}
\usepackage{amsthm,amsmath}
\usepackage{latexsym}
\makeatletter
\renewcommand{\@seccntformat}[1]
{\csname the#1\endcsname.\enspace}
\makeatother
\begin{document}
\newtheorem{theorem}{Theorem}
\newtheorem{lemma}{Lemma}
\newtheorem{remark}{Remark}
\newtheorem{proposition}{Proposition}
\newtheorem{corollary}{Corollary}
\newtheorem{example}{Example}
\newtheorem{definition}{Definition}

\begin{center}
\textbf{Estimation of a nonnegative location parameter with unknown scale} 
\\
\end{center}

\begin{center}
{Mohammad Jafari Jozani$^{a,}$\footnote{Corresponding author: m$_{-}$jafari$_{-}$jozani@umanitoba.ca}, \' Eric
 Marchand$^{b}$, and William E. Strawderman$^{c}$ } \\
\end{center}
\begin{center}
{\it a University of Manitoba, Department of Statistics, Winnipeg, MB, CANADA, R3T 2N2 } \\
{\it b  Universit\'e de
    Sherbrooke, D\'epartement de math\'ematiques, Sherbrooke, QC,
    CANADA, J1K 2R1} \\
    {\it c Department of Statistics and Biostatistics, Rutgers University, 561 Hill
Center,
 Busch Campus, Piscataway, N.J., USA 08854-8019 
} \\
\end{center}

\begin{center}
{\sc Summary}
\end{center}
\small  
For normal canonical models, and more generally a vast array of general spherically symmetric location-scale models with a residual vector, we consider estimating
the (univariate) location parameter when it is lower bounded.  We provide conditions for estimators to dominate the benchmark minimax MRE estimator, and thus be minimax under scale invariant loss.  These minimax estimators include the generalized Bayes estimator with respect to the truncation of the common non-informative prior onto the restricted parameter space for normal models under general convex symmetric loss, as well as non-normal models under scale invariant $L^p$ loss with $p>0$.  We cover many other situations when the loss is asymmetric, and where other generalized Bayes estimators, obtained with different powers of the scale parameter in the prior measure, are proven to be minimax.  We rely on various novel representations, sharp sign change analyses, as well as capitalize on Kubokawa's integral expression for risk difference technique.  Several other analytical properties are obtained, including a robustness property of the generalized Bayes estimators above when the loss is either scale invariant $L^p$ or asymmetrized versions.  Applications include inference in two-sample normal model with order constraints on the means.


\normalsize

\section{Introduction}

\subsection{Preamble}
We begin with the normal model in canonical form  
\begin{equation}\label{model}
X\sim~N(\mu,\sigma^{2}),S^{2}\sim\sigma^{2}\chi^{2}_{n},\quad\textrm{independent}\quad
(n\geq1),
\end{equation}
which plays a central role in both statistical theory and practice.  Consider situations where additional
information on $(\mu,\sigma)$ is available in terms of parametric restrictions.  Bayesian inference in such restricted
parameter space problems does not, conceptually, present any difficulties as both the prior and the resulting posterior will be
adapted and will adapt to the constraints.  Assessing the frequentist performance of Bayesian estimators in such situations is, however,
considerably more challenging.  Such assessments may include, for instance, testing for minimaxity, an evaluation in comparison
to a benchmark procedure such as minimum risk equivariant (MRE) estimator or a maximum likelihood estimator (mle), or a study of the frequentist performance
of associated Bayesian confidence intervals.
 
As an illustration, consider model (\ref{model}) with known $\sigma$ and the nonnegative mean restriction $\mu \geq 0$. 
Despite early discoveries by Katz (1961) and Sacks (1963) that the generalized Bayes estimator with respect to the flat prior
on $[0,\infty)$ is minimax and dominates the MRE estimator $\delta_0(X,S)=X$ under squared error loss, despite various generalizations to other models and location invariant losses (Farrell, 1964; Kubokawa, 2004; Marchand and Strawderman, 2005), no {\it other} Bayes minimax estimators were known until the Maruyama and Iwasaki (2005) findings which provide other Bayes minimax estimators under squared error loss.  Even then, little has been obtained for estimating $\mu$ in (\ref{model}) for $\mu \geq 0$ and unknown $\sigma$.  In this case, Kubokawa (2004) obtained, for scale invariant squared error loss, a class of minimax improvements on $\delta_0$, which includes the generalized Bayes estimator $\delta_{\pi_0}(X,S)$ with respect the truncation of the usual non-informative prior onto the restricted parameter space (see expression \ref{prior}).

Our main motivation for his work has been to generalize and better understand Kubokawa's findings.  The paper consists of various extensions with respect to the loss, the model, and the prior; which bypass in a unified way the specific normal case-squared error loss calculations by Kubokawa.  Several new technical aspects have been developed to meet such challenges.

\subsection{The problem}

As an extension of model (\ref{model}), we consider spherically symmetric models for an observable $(X,U)=(X,U_1, \ldots, U_n)$ with density proportional to
\begin{equation}
\label{XU}  \frac{1}{\sigma^{n+1}} \, f(\frac{(x-\mu)^2 + \|u\|^2}{\sigma^2})\,,
\end{equation}
and with $n \geq 1$, $\mu \geq 0$, $\sigma >0$.  The function $f: \mathbb{R}^+ \to \mathbb{R}^+$ is known, and it assumed throughout that:
\begin{equation}
\label{assumptions}
f'<0, \; \hbox{and} \;\, \frac{t f'(t)}{\;f(t)} \;\, \hbox{decreases in} \; t \; \hbox{for} \; t>0\,.  
\end{equation}
Hereafter, for conciseness, reference to model (\ref{XU}) shall be understood to encompass these assumptions on $f$.
Multivariate (for $X$) versions of (\ref{XU}) have been previously considered, namely in recent work where robust minimax generalized Bayes estimators of $\mu$ without constraints are provided (see Fourdrinier and Strawderman, 2010).  Various other features of such models are described in Section 2.1.  

We consider estimating $\mu$ where it is assumed that $(\mu, \sigma) \in \Theta=\{(\mu,\sigma): \mu \geq 0, \sigma >0\}$
under location and scale invariant loss
\begin{equation}\label{loss}
\rho(\frac{d-\mu}{\sigma}),
\end{equation}
with (i) $\rho$ absolutely continuous a.e., (ii) $\rho$ strictly bowled shaped with
$\rho(t)\geq\rho(0)=0$ for all $t \in \mathbb{R}$, $\rho'<0$ on $(-\infty,0)$ and $\rho'>0$ on $(0,\infty)$.  We also assume that the pair $(f,\rho)$ leads to risk finiteness,
namely that there exists a unique minimum risk equivariant estimator for the unconstrained problem.  In such cases, it is given by $\delta_{0}(X,S)=X+c_{0}S$ with
constant risk $R((\mu,\sigma), \delta_0)= E_{0,1}[\rho(X+c_{0}S)]$, and with (also see Remark \ref{mreremark})
\begin{equation}
\label{mre}  c_{0}=\hbox{argmin}_c \{E_{0,1}[\rho(X+cS)]\}\,,
\end{equation}
which is uniquely determined by 
$E_{0,1}[S\rho'(X+c_{0}S)]=0 \,.$
It is also worth pointing out that $c_{0}=0$ for symmetric losses $\rho$, and consequently that the MRE estimator coincides with the unbiased estimator $X$, and is robust with respect to the choice of the underlying model density $f$.
It follows from Kiefer (1957) that $\delta_0$ is minimax for the unconstrained problem.  With the constraint on $\mu$, $\delta_0(X,S)$ produces indeed 
implausible estimates, but it remains minimax (see Marchand and Strawderman, 2011, and references therein) for general $\rho$, and its constant risk thus matches the minimax risk.   The challenge here is to search for good improvements on $\delta_{0}(X,S)$ that capitalize
on the parametric information, and we focus on potential Bayesian improvements such as the generalized
Bayes estimators $\delta_{\pi_l}$ with respect to the prior measures
\begin{eqnarray}\label{prior}
\pi_{l}(\mu,\sigma)=\frac{1}{\sigma^{l+1}} \, 1_{[0, \infty)} (\mu)1_{(0,\infty)}(\sigma)\,; l \geq -(n-1)\,;
\end{eqnarray}
the lower bound on $l$ required for the posterior density to be well defined.
The class includes the choice $\pi_0$ which is of intrinsic interest as it represents a plausible adaptation,
or truncation onto $\Theta$ of the right Haar invariant measure
$\pi_{rh}$ with the MRE estimator (also) being the generalized Bayes
estimator $\delta_{\pi_{rh}}$ with respect to $\pi_{rh}$.  Moreover, the study of frequentist properties on the restricted parameter space of Bayesian procedures associated with $\pi_0$ or, more generally, truncations of the right Haar invariant prior measure has recently surfaced in interval estimation problems (Zhang and Woodroofe, 2003; Marchand and Strawderman, 2006; Marchand et al., 2008).

In Section 2, we further describe features of the underlying model and present various expressions, properties, and illustrations relative to the Bayes estimators $\delta_{\pi_l}$.  Namely, we establish a robustness property, applicable to scale invariant $L^p$ loss with $\rho(t)=|t|^p$, $p>0$, and asymmetrized versions as given in (\ref{lp}), stating that the Bayes estimator $\delta_{\pi_l}$ does not depend on the underlying $f$ in (\ref{XU}).

The developments of Section 3 make use of Kubokawa's (1994) IERD (Integral Expression
of Risk Difference) technique to derive classes of dominating (minimax) estimators of $\delta_0(X,S)=X+c_0S$.  With further analyses, which bring into play novel
technical arguments of interest on their own, we provide several instances of
$(f,\rho)$ where these classes of  minimax estimators include Bayesian estimators of the type $\delta_{\pi_l}$.  Namely, we establish in Sections 4 and 5 that:

\begin{enumerate}
\item[ {\bf (A)}]  The Bayes estimators $\delta_{\pi_l}$ with $l \geq 0$ dominate $\delta_0$ for normal models in (\ref{model}) and general convex $\rho$'s such that
$\rho$ is even.  The estimator $\delta_{\pi_0}$ also dominates $\delta_0$ for asymmetric $\rho's$ such that $|\rho'(u)| \geq |\rho'(-u)|$ for all $u>0$;

\item[ {\bf (B)}]  The Bayes estimators $\delta_{\pi_l}$ with $l \geq 0$ dominate $\delta_0$ for all (fixed) $f$ in (\ref{XU}) satisfying assumption (\ref{assumptions}), and whenever the loss is scale invariant $L^p$, $p \geq 1$.   The estimator $\delta_{\pi_0}$ also dominates $\delta_0$ for asymmetrized versions as given in (\ref{lp}) (where $|\rho'(u)| \geq |\rho(-u)|$ for all $u>0$ as in {\bf (A)};

\item[ {\bf (C)}]  The Bayes estimator $\delta_{\pi_0}$ dominates $\delta_0$ for all (fixed) $f$ in (\ref{XU}) satisfying assumption (\ref{assumptions}), and whenever the loss is scale invariant $L^p$ with $p \in (0,1)$.

\end{enumerate}

The ensemble of results provide extensions of Kubokawa's normal case, scale invariant squared error loss result applicable to $\delta_{\pi_0}$ in three directions:
choice of $f$, choice of $\rho$, and applicability to other Bayesian estimators $\delta_{\pi_l}$'s.  Moreover, the developments relative to {\bf (A), (B)}, and {\bf (C)} 
are unified and contain two alternative proofs replicating Kubokawa's result.  It is also notable that {\bf (C)} involves the case of a concave in $|\frac{d-\mu}{\sigma}|$ (and hence non-convex) loss. 
Finally, various other observations, including non-minimaxity results, are also given throughout the exposition and in Section 6.

\section{Preliminary results and properties of the estimator $\delta_{\pi_l}$}

\subsection{The underlying model}  

In (\ref{XU}) and (\ref{assumptions}), the density of $(X,U)$ is unimodal with central location parameter 
$(\mu, 0, \ldots, 0)$ and scale parameter $\sigma$.  Our parameter of interest is the nonnegative $\mu$, or median, of $X$, while $U$ is a residual vector.
Condition (\ref{assumptions}) is equivalent to an increasing monotone likelihood ratio (mlr) in $(X-\mu)^2 + \|U\|^2$ of the family of densities in
(\ref{XU}) when viewed as a scale family (parameter $\sigma$) with known $\mu$.  Assumption (\ref{assumptions}) is, for unimodal and symmetric densities, weaker than
both {\bf (a)} the logconcavity of $f(y)$ and {\bf (b)} the logconcavity of $f(y^2)$ for $y>0$, with {\bf (a)} implying {\bf (b)}, and with {\bf (b)} equivalent to an increasing mlr property in $X$ of the family of densities in (\ref{XU}) when viewed as a location family (parameter $\mu$) with known $\sigma$.  

The most important and best known case covered by (\ref{XU}) and (\ref{assumptions}) is the normal case where   
$(X,U) \sim N_{n+1}((\mu,0, \ldots, 0), \sigma^2 I_{n+1})$ and $f(t) \propto \, e^{-t/2}$.  However, our inference results will also apply to many other 
models such as (i) exponential power densities with $f(t) \propto e^{-\alpha t^p}$, $p>0$, $\alpha>0$, including Laplace densities arising for $p=1/2$; (ii) the Kotz distribution with $f(t) \propto t^m  e^{-\alpha t}$, $m \in (-1/2,0)$, $\alpha>0$; as well as for (iii) Student densities with $f(t) \propto (1+t/\nu)^{-(\nu +n+1)/2}$, $\nu \geq 1$ degrees of freedom.   The Student example illustrates a non-logconcave $f$ (in fact, it is logconvex) which satisfies the weaker assumptions required here.  The Student distributions, which are scale mixtures of normals, often serve as useful, alternative models to the normal model.  Here is an interesting general situation for which scale mixtures inherit assumption (\ref{assumptions}). 

\begin{lemma}
\label{mixtures}
A scale mixture of the form $f(t)= \int_0^{\infty} v f_0(tv) h(v) dv$ satisfies assumption (\ref{assumptions}) as soon as both $f =f_0$ and $f=h$ satisfy assumption (\ref{assumptions}).
\end{lemma}
{\bf Proof.}  See Appendix.

\begin{remark}
In the Student case above, both $f_0$ (a normal density) and $h$ (a gamma density), are logconcave and satisfy (\ref{assumptions}).
\end{remark}

\begin{remark}
We note that model (\ref{XU}) arises for observables $Y_1, \ldots, Y_{n+1}$ having joint density $$\frac{1}{\sigma^{n+1}}
f\left(\frac{\sum_i(y_i-\theta)^2}{\sigma^2}\right)
,$$
through an orthogonal transformation $$(Y_1, \ldots, Y_{n+1}) \to (X=\sqrt{n} \, \bar{Y}, U_1, \ldots, U_n),$$ with $\mu=\sqrt{n} \, \theta.$
\end{remark}

For model (\ref{XU}), $(X, S=\|U\|)$ is a sufficient statistic with joint density $f_{X,S}$ on
$\mathbb{R} \times \mathbb{R}^+$ which we take as equal to:
\begin{equation}
\label{reduced}
\frac{s^{n-1}}{\sigma^{n+1}} \, f(\frac{(x-\mu)^2 + s^2}{\sigma^2})\,.
\end{equation}

For the normal model canonical form in (\ref{model}),  
we will write the joint density of $(X,S)$ in (\ref{reduced}) as
$\frac{1}{\sigma^2} \phi(\frac{x-\mu}{\sigma}) \, h(\frac{s}{\sigma})$, with
\begin{equation}
\label{densities}
\phi(u)=(2\pi)^{-1/2} e^{-u^2/2}\,, \; \hbox{and} \; h(v)=\frac{v^{n-1}\, e^{-v^2/2}}{\Gamma(n/2) \,2^{n/2-1}}\,.
\end{equation}

\subsection{Properties of the Bayes estimators $\delta_{\pi_l}$}

We proceed with various preliminary results, observations, and illustrations concerning the generalized Bayes estimators $\delta_{\pi_l}(X,S)$ with respect to the improper priors in (\ref{prior}).  As previously mentioned, one can verify that the lower bound on the power $l$ in (\ref{prior}) guarantees that the posterior density of $(\mu,\sigma)$ is well defined given that
(\ref{reduced}) is a density.  Even with a well defined posterior density, we further assume, and not 
necessarly emphasize (mainly in Sections 3,4, and 5), that the pair $(f,\rho)$ leads to the existence of the Bayes estimator $\delta_{\pi_l}$.  

We define for $m>0, w \in \mathbb{R},z \in \mathbb{R}$,
\begin{equation}
\label{B}
B_m(w,z)= \int_0^{\infty} \int_{-\infty}^{vw} \rho'(u+c_0v + zv) \,v^{m} \, f(u^2+v^2) \, du \, dv\,,
\end{equation}
provided it exists.  The function $B_m(w,z)$, as well as some of its properties (see for instance
Lemma \ref{Blemma}) will play a key role below, namely in the following
representation of the Bayes estimator $\delta_{\pi_l}(X,S)$.

\begin{lemma}
\label{deltapi*}  Under model (\ref{reduced}), provided existence of the Bayes estimator $\delta_{\pi_l}$, we have $\delta_{\pi_l}(X,S)=X+ c_0S + g_{\pi_l}(\frac{X}{S})S $, where
$g_{\pi_l}(y)$ satisfies, for all $y \in \mathbb{R}$, $l > -(n-1)$,
\begin{equation}
\label{deltapi*def}
B_{n+l}(y, g_{\pi_l}(y))=0.
\end{equation}
\end{lemma}
{\bf Proof.}
Writing an estimator as $X+ c_0 S +g(X,S)$, we have that the Bayes estimate
$\delta_{\pi_l}(x,s)$ minimizes in $g(x,s)$ the expected posterior loss:
$$ E[\rho(\frac{x+c_0s+ g(x,s) -\mu }{\sigma})|(X,S)=(x,s)], $$
or, equivalently,
$$ \int_{0}^{\infty} \int_{0}^{\infty} \rho(\frac{x+c_0s+ g(x,s) -\mu
}{\sigma})
\frac{s^{n-1}}{\sigma^{n+1}} f(\frac{(x-\mu)^2 + s^2}{\sigma^2})
\frac{1}{\sigma^{l+1}} d\mu d\sigma\;. $$
With the change of variables $(\mu, \sigma) \rightarrow
(u=\frac{(x-\mu)}{\sigma},
v=\frac{s}{\sigma})$, the Bayes
estimate
$\delta_{\pi_l}(x,s)$ is seen to minimize in $g(x,s)$:
$$\int_{0}^{\infty} \int_{-\infty}^{v x/s} \rho(u+c_0v+ \frac{v}{s} g(x,s))
f(u^2+v^2) \, v^{n+l-1} \, du dv\;. $$
Now, observe that $\frac{1}{s} g(x,s)$ depends on $(x,s)$ only through the
function $y=x/s$, which implies that the estimator $\delta_{\pi_l}(X,S)$
is of the form $X+c_0S + g_{\pi_l}(\frac{X}{S})\,S$ with $g_{\pi_l}(y)$
minimizing in $g(y)$ the quantity
\begin{equation}
\label{delta0def}
\int_0^{\infty} \, \int_{-\infty}^{vy} \rho(u+c_0v + g(y)v) \, f(u^2+v^2) \, v^{n+l-1} \,du dv\,.
\end{equation}
Finally, the result is obtained by differentiation.  \qed

We point out that $g_{\pi_l}(y)$ is uniquely determined (Lemma \ref{Blemma}), and is a continuous
function of $y$ such that
\begin{equation}
\label{positive}
g_{\pi_l}(y) \geq -y -c_0 \;\; \hbox{for all} \; y \in \mathbb{R}.
\end{equation}
This must indeed be the case as the Bayes estimates
$\delta_{\pi_l}(x,s)$ are necessarily nonnegative, and
$\frac{1}{s} \delta_{\pi_l}(x,s) \geq 0
\Longleftrightarrow
\frac{x}{s} +c_0 + g_{\pi_l}(\frac{x}{s}) \geq 0 $. 
We pursue with an intriguing robustness property, and alternative representation, of the Bayes
estimators $\delta_{\pi_l}$ for scale invariant $L^p$ loss, and their asymmetrized versions given by
\begin{equation}
\label{lp}\rho_{c_1,c_2}(t)= c_1 \,|t|^p \,1_{(-\infty,0)}(t) \,+ \,c_2 \,|t|^p \,1_{[0,\infty)}(t)\,, 
\end{equation}
with $p>0$, $c_1>0$, and $c_2>0$.
\begin{lemma}
\label{independence}
For losses $\rho_{c_1,c_2}$ as in (\ref{lp}), the Bayes estimators $\delta_{\pi_l}$, given in Lemma 
\ref{deltapi*}, do not depend on the underlying model density $f$ provided they exist.
\end{lemma}
{\bf Proof.}  
From (\ref{delta0def}), we have
\begin{align}\label{argmin}
c_0+ g_{\pi_l}(y)
&=  \hbox{argmin}_h \int_0^{\infty} \int_{-\infty}^{Vy} \, \rho_{c_1,c_2}(u+hv) \, 
f(u^2+v^2) \, v^{n+l-1} \, du dv  \nonumber\\ 
 &= \hbox{argmin}_h \int_0^{\infty} \int_{-\infty}^{Vy} \, \rho_{c_1,c_2}(\frac{u}{v}+h) \, 
f(u^2+v^2) \, v^{n+l+p-1} \, du dv\nonumber \\
 &= \hbox{argmin}_h \left( \int_0^{\infty} x^{(n+l+p-1)/2} f(x) \, dx \right) \, \left( \int_{-\infty}^{y} \, 
\frac{\rho_{c_1,c_2}(t+h)}{(1+t^2)^{(n+l+p+1)/2}} \, dt \right)  \\  
\label{argmin} &= \hbox{argmin}_h \,  \, \int_{-\infty}^{y} \, 
\frac{\rho_{c_1,c_2}(t+h)}{(1+t^2)^{(n+l+p+1)/2}} \, dt \, ,  \nonumber
\end{align}
by making use of the homogeneity of $\rho_{c_1,c_2}$ and the change of variables $(u,v) \to (t=u/v, x=u^2+v^2)$.
Finally, expression (\ref{argmin}) tells us that $\delta_{\pi_l}(x,s)=x+s(c_0+g_{\pi_l}(x/s))$ is independent of $f$.  \qed\\

This type of property seems to have first been noticed by Maruyama (see Maruyama, 2003; Maruyama and Strawderman, 2005) in a multivariate setting under $L^2$ loss.

\begin{remark} (Minimum risk equivariant estimator)
\label{mreremark}
\begin{enumerate}
\item[ {\bf (a)}]
Proceeding as in the proof of Lemma \ref{deltapi*}, we obtain the useful representation $X+ c_0(n) S$ for the MRE estimator, with the defining equation
\begin{equation}
\label{c0(m)}
\int_0^{\infty} \int_{-\infty}^{\infty} \, \rho'(u+c_0(m) v) \, v^m \, f(u^2+v^2) \, du \,dv\,=\,0\,,
\end{equation}
for $c_0(m)$, $m \geq 1$.  

\item[ {\bf (b)}] A robustness property similar to Lemma \ref{independence}
(also illustrated in Example \ref{L1L2}, part ({\bf C}) is shared by the MRE estimators with respect
to losses $\rho_{c_1,c_2}$ and can be established
by expanding (\ref{mre}) showing that 
\begin{equation}
\label{mrelp}
c_0=\hbox{argmin}_c \,  \, \int_{-\infty}^{\infty} \, 
\frac{\rho_{c_1,c_2}(t+c)}{(1+t^2)^{(n+p+1)/2}} \, dt \,.
\end{equation}
\end{enumerate}
\end{remark}

\begin{example} (scale invariant $L^2$ loss, scale invariant $L^1$ loss and their asymmetrized versions)
\label{L1L2}
\begin{enumerate}
\item[ \bf{(A)}]  For scale invariant squared error loss with $\rho(t)=t^2$ in (\ref{loss}), the MRE 
estimator is $\delta_0(X)=X$, provided the second moment of $X$ under (\ref{XU}) exists.  Lemma \ref{deltapi*} as well as (\ref{argmin}) provide representations $X+g_{\pi_l}(\frac{X}{S})S$ for the 
Bayes estimator $\delta_{\pi_l}(X,S)$; $l > -(n-1)$.  Differentiating (\ref{argmin}) with respect to $h$, we obtain directly for $y \in \mathbb{R}$
\begin{equation}
\label{gpill}
g_{\pi_l}(y)= - E[T|T \leq y ]\,,
\end{equation}
where $T$ has density on $\mathbb{R}$ proportional to $(1+t^2)^{-(n+l+3)/2}$.  Here the distribution of $T$ is a multiple of a Student distribution with $n+l+1$ degrees of freedom.  Equivalently from (\ref{deltapi*def}), we have 
\begin{eqnarray}
\nonumber \, & & B_{n+l}(y,g_{\pi_l}(y))=0 \\
\nonumber \Longleftrightarrow & &\int_0^{\infty} \int_{-\infty}^{vy} (u + g_{\pi_l}(y) v) \, v^{n+l} \, f(u^2+v^2) \,du\,dv =0  \\
\label{gpil}\Longleftrightarrow & & g_{\pi_l}(y)= -\frac{\int_0^{\infty} \int_{-\infty}^{vy} \frac{u}{v} \, v^{n+l+1} \, f(u^2+v^2) \,du\,dv}
{\int_0^{\infty} \int_{-\infty}^{vy} \, v^{n+l+1} \, f(u^2+v^2) \,du\,dv}\,, 
\end{eqnarray}
illustrating the fact that the distribution of $T$ arises as the (independent of $f$) distribution of the ratio $\frac{U}{V}$, with $(U,V)$ having joint density on $\mathbb{R} \times \mathbb{R}^+$ proportional to $v^{n+l+1} f(u^2+v^2)$. 
From representation (\ref{gpill}), observe that $g_{\pi_l}(\cdot)$ decreases on $\mathbb{R}$ with $\lim_{y \to \infty} g_{\pi_l}(y)=0$ 
(since $\int_{-\infty}^{\infty} u f(u^2+v^2) du =0$ for all $v>0$), and hence that $g_{\pi_l}(\cdot)$ is positive, i.e., 
$\delta_{\pi_l}$ expands on the MRE $\delta_0$.
Such properties are of interest as they indicate that the amplitude of the expansion $\delta_{\pi_l}(x,s) - \delta_0(x,s)$ decreases in $x$ for fixed $s$, and increases in $s$ for fixed $x$ (in fact $(\delta_{\pi_l}(x,s) - \delta_0(x,s))/s$ increases in $s$).  Such a property resonates back to Katz (1961) where in the normal case with known $\sigma$, the Bayes estimator with respect to a flat prior for $\mu$ on $(0,\infty)$ expands $X$ by the amount $\sigma \frac{\phi(x/\sigma)}{\Phi(x/\sigma)}$ which decreases in $x$ and increases in $\sigma$.   Below, we establish such properties for general convex $\rho$ in Lemma \ref{gpi*}, as well as scale invariant $L^p$ concave loss with $p \in (0,1)$ in Lemma \ref{gpi*}.   
Finally, we point out that alternative expressions for $\delta_{\pi_1}$ in the above normal case were given by Kubokawa (2004), as well as Marchand, Jafari Jozani, and Tripathi (2011). 

\item[ \bf{(B)}]  As above, for scale invariant absolute value error loss with $\rho(t)=|t|$ in (\ref{loss}), the MRE estimator is $\delta_0(X)=X$.
For $l \geq -(n-1)$, $\delta_{\pi_l}(X,S)=X+g_{\pi_l}(\frac{X}{S})S$ is obtainable from (\ref{argmin}) yielding 
\begin{equation}
\label{gpiL1} g_{\pi_l}(y)= - \hbox{median}[T|T \leq y ] =  -F_{n+l}^{-1}(\frac{F_{n+l}(y)}{2})\,,
\end{equation}
where $F_m$ and $F_m^{-1}$ are the cdf and inverse cdf of $T$ having density on $\mathbb{R}$ proportional to $(1+t^2)^{-(m+2)/2}$. 
As above, it is easily seen directly that such a $g_{\pi_l}(\cdot)$ decreases on $\mathbb{R}$, that $\lim_{y \to \infty} g_{\pi_l}(y)=0$, 
that $\delta_{\pi_l}(x,s)$ expands once again on $\delta_0(x,s)$ for all $(x,s) \in \mathbb{R} \times \mathbb{R}^+$, and the difference between these estimates decreases in $x/s$.

\item[ \bf{(C)}]  Consider now asymmetrized $L^1$ losses $\rho_{c_1,c_2}$ in (\ref{lp}) with $p=1$.  By making use of Remark \ref{mreremark}, the MRE estimator is given by $\delta_0(X)=X+c_0S$, with $c_0$ independent of $f$, and $c_0(n)=-F_n^{-1}(\frac{c_2}{c_1+c_2})$ and $F_n^{-1}$ the inverse cdf given in part ({\bf B}).
For $l \geq -(n-1)$, we obtain from (\ref{argmin}) $\delta_{\pi_l}(X,S)=X+c_0(n)S+ g_{\pi_l}(\frac{X}{S})S$  
with $g_{\pi_l}(y)=-c_0(n) -F_{n+l}^{-1}(\frac{c_2}{c_1+c_2} F_{n+l}(y))$, thus extending (\ref{gpiL1}) which occurs for $c_1=c_2$.
Observe here that $\lim_{y \to \infty} g_{\pi_l}(y)=-c_0(n) + c_0(n+l)$, which does not equal $0$ in general, the exception being precisely $l=0$, and/or $c_1=c_2$.
This property is more general as seen below in Lemma \ref{gpi*}.
\end{enumerate}
\end{example}

We pursue with further properties relative to $B_m(\cdot,\cdot)$ and $g_{\pi_l}$ (applicable when these quantities exist).

\begin{lemma}
\label{Blemma}
For all $a > 0$, $y \in \mathbb{R}$, $l \geq -(n-1)$, and strictly bowled-shaped $\rho$,

\begin{enumerate}

\item[ {\bf (a)}]

$ B_{n+l}(y+a, g_{\pi_l}(y)) >0$;

\item[{ \bf (b)}]  $B_{n+l}(y,z)$ is nondecreasing in $z$ whenever $\rho$ is also convex;

\item[ {\bf (c)}]  $\lim_{y \rightarrow \infty} B_{n+l}(y,0)=0$ whenever $l=0$; or whenever $l \neq 0$ and $\rho$ is an even function.

\end{enumerate}
\end{lemma}

{\bf Proof.}
Part (b) is obvious given the convexity of $\rho$, while part (c)
follows from the given representations (\ref{c0(m)}) and (\ref{deltapi*def}).
For establishing (a), suppose, in order to arrive at a contradiction, that
$B_{n+l}(y+a, g_{\pi_l}(y)) \leq 0. $  This would imply $C_1 \leq 0$,
where
$$ C_1=\int_0^{\infty}  \int_{vy}^{v(y+a)} \rho'(u+c_0v+g_{\pi_l}(y)v) \; v^{n+l} f(u^2+v^2) \; du \, dv \,.$$
Now, observe that for $(u,v) \in I(u,v)=\{ (u,v): vy < u <v(y+a)  \}$,
we have by
(\ref{positive}): $u+c_0v +g_{\pi_l}(y)v > vy +c_0v+g_{\pi_l}(y)v \geq 0$,
implying $\rho'(u+ c_0v +g_{\pi_l}(y) v) > 0$, (for such $(u,v)'s \in
I(u,v)$).
This renders $C_1 \leq 0$ impossible, and yields the result. \qed

The strictly decreasing property of $g_{\pi_l}$ that follows in Lemma \ref{gpi*} is a critical property that we will exploit later for the risk comparisons.  We do not know how far the property can be extended for non-convex $\rho$, but we do establish here, and use later, such a property for $L_p$ losses and their asymmetrized versions for the non-convex choices $p \in (0,1)$.

\begin{lemma}
\label{gpi*}
For $l \geq -(n-1)$,
\begin{enumerate}
\item[ {\bf (a)}]  $g_{\pi_l}(y)$ is strictly decreasing in $y$ whenever $\rho$ is convex;
\item[ {\bf (b)}]   $g_{\pi_l}(y)$ is strictly decreasing in $y$ whenever the loss is
$\rho_{c_1,c_2}$ as in (\ref{lp}) with $p \in (0,1)$. 
\item[ {\bf (c)}]  For strictly bowled-shaped $\rho$, $\lim_{y \rightarrow \infty} g_{\pi_l}(y)=-c_0(n)+ c_0(n+l)$, where $c_0(m)$ is defined 
in (\ref{c0(m)}).  Consequently, $\lim_{y \rightarrow \infty} g_{\pi_l}(y)=0$ whenever $l=0$, or $l \neq 0$
and $\rho$ is even.
\end{enumerate}
\end{lemma}
{\bf Proof.}
{\bf (a)}
It suffices to show that we cannot have
$g_{\pi_l}(y+\epsilon) \geq g_{\pi_l}(y)$ for some
$y \in \mathbb{R}$, $\epsilon >0$.  Indeed, if this were the case, it would follow,
using defining equation (\ref{deltapi*def}) and part (a) of Lemma \ref{Blemma}, that
$$  0 = B_{n+l}(y+\epsilon, g_{\pi_l}(y+\epsilon)) \geq
B_{n+l}(y+\epsilon, g_{\pi_l}(y)) >0 \;,$$ which is not possible.

{\bf (b)} Set $s(y)=-c_0-g_{\pi_l}(y)$ and rewrite 
representation (\ref{argmin}) as
\begin{equation}
\label{s}
s(y)= \hbox{argmin}_s E[\rho_{c_1,c_2}(T-s)|T \leq y],
\end{equation}
with $T$ having density proportional to $(1+t^2)^{-(n+l+p+1)/2}$ on $\mathbb{R}$.  Observe that the family of densities
for $T|T \leq y$ has strictly increasing monotone likelihood ratio in $T$ with parameter $y$.  Now, consider, for $a_1 <a_2$, the 
function $\rho_{c_1,c_2}(t-a_1) - \rho_{c_1,c_2}(t-a_2)$, which changes signs once from $-$ to $+$ as a function of $t$ as $t$ increases on 
$\mathbb{R}$, and infer that
$$ H(a_1,a_2,y)=E[\rho_{c_1,c_2}(T-a_1) - \rho_{c_1,c_2}(T-a_2)] $$
has a single root, and changes signs once from $-$ to $+$, as a function of $y$, as $y$ increases on 
$\mathbb{R}$, given the mlr property (e.g., Lehmann, 1986).  
Suppose now, in order to arrive at a contradiction that $g_{\pi_l}$ is not strictly decreasing, i.e., $s$ is not strictly increasing  
and there exists $y_2<y_1$ such that $a_2=s(y_2) \geq s(y_1)=a_1$.  Then, we would have with the definition of $s(y)$ in (\ref{s}) and the properties of
$H$: $H(s(y_1), s(y_2),y_2) >0$ and $H(s(y_1), s(y_2), y_1) <0$ which leads to a contradiction and establishes the result.

{\bf (c)}  This follows by matching expression (\ref{deltapi*def}) when $y \to \infty$ with 
(\ref{c0(m)}).  \qed

\begin{remark}
The above proof in ({\bf b}) goes through for all losses $\rho_{c_1,c_2}$, including the convex cases with $p \geq 1$.
\end{remark}

The following results permit the ordering of Bayes estimators $\delta_{\pi_l}$ in terms of the power $l$ in the prior measure
$\pi_l$ in (\ref{prior}).

\begin{lemma}
\label{l}  For the normal model in (\ref{model}), $y \in \mathbb{R}$, and convex and even $\rho$, the quantities $g_{\pi_l}(y)$ decreases in $l$, $l \geq -(n-1)$, provided they exist.
\end{lemma}
{\bf Proof.}  See Appendix.  
\begin{corollary}
\label{c}
For models (\ref{XU}) with $f$ satisfying assumption (\ref{assumptions}), $y \in \mathbb{R}$, and scale invariant $L^p$ loss with $p>0$, $g_{\pi_l}(y)$ decreases in $l$, $l \geq -(n-1)$, provided existence. 
\end{corollary}
{\bf Proof.}  Lemma \ref{independence} tells us that $g_{\pi_l}(y)$ is independent of $f$ and thus matches the normal model $g_{\pi_l}(y)$ and
Lemma \ref{l} tells us that such $g_{\pi_l}(y)$'s decrease in $l$ whenever $\rho$ is even as for the $L^p$ loss here.  \qed

\section{Minimax Conditions for general $\rho$ and $f$}

For estimating $\mu \geq 0$ in (\ref{XU}) or in (\ref{reduced}) with unknown $\sigma>0$ under strictly bowled-shaped loss $\rho(\frac{d-\mu}{\sigma})$, 
we establish here useful sufficient conditions for an estimator $\delta(X,S)$ to be minimax.  We first make use of Kubokawa's IERD technique in Theorem
\ref{main}.  Proposition \ref{main2} (below) then extracts a sign varying condition for minimaxity which will serve as the basis for further analysis for
the specific cases of normal models and general convex $\rho$ in Section 4, and for $L^p$ losses and their asymmetric versions $\rho_{c_1,c_2}$ with general $f$
satisfying assumptions (\ref{assumptions}) in Section 5.  Various other technical results and remarks, including a condition for non-minimaxity with applications, are
also introduced in this section. 
We consider the following subclass of scale invariant estimators.
\begin{definition}
\label{class}
$C=\{\delta_g(X,S): \delta_g(X,S)=\delta_0(X,S) + g(\frac{X}{S}) \,S $,
with $g$ absolutely continuous a.e., nonincreasing, non-constant, and 
$\lim_{t \rightarrow \infty} g(t)=0 \}$.
\end{definition}
These estimators in $C$ expand upon $\delta_0$, in view of the restriction $\mu \geq 0$, include $\delta_{\pi_0}$ and the
generalized Bayes estimators $\delta_{\pi_l}$; $l \neq 0$, $l \geq -(n-1)$; for even $\rho$ as seen by the properties given in Lemma \ref{gpi*}.  Under invariant losses as in (\ref{loss}), such estimators will have frequentist risk $R(\theta,\delta_g)$ depending on $\theta=(\mu,\sigma)$ only through the maximal invariant $\lambda=\mu/\sigma$, and we seek conditions for which such a risk falls below the constant risk of the MRE estimator $\delta_0$ for all $\lambda \geq 0$.  As mentioned above, such improvements will necessarily be  minimax estimators since $\delta_0$ is minimax.  Hereafter, we will just refer, for the most part, to such improvements as being minimax estimators.  
The focus is largely on the generalized Bayes estimator $\delta_{\pi_0}$, which will be seen to be minimax for various settings of $(f,\rho)$ and which provides a benchmark
in the sense that estimators $\delta_g \in C$ will be minimax for convex $\rho$ under the simple condition that $\delta_g$ not expand on $\delta_0$ as much as $\delta_{\pi_0}$ (Theorem \ref{main}, {\bf (ii)}).  In turn, for various choices of $(f,\rho)$ with $\rho$ even, and by appealing to Lemma \ref{l}, these classes of minimax estimators will contain the generalized Bayes estimators $\delta_{\pi_l}$'s, $l>0$.  We now pursue with an intermediate dominance condition. 

\begin{theorem}
\label{main}
For estimating $\mu$ in (\ref{XU}) or (\ref{reduced}) with $\mu \geq 0, \sigma >0$, an estimator
$\delta_g \in C$ is minimax, under
strictly bowled shaped loss $\rho(\frac{d-\mu}{\sigma})$ whenever either one of the following conditions
holds for all $\lambda \geq 0$ and $y \in \{y:g'(y) < 0\}$: 
\begin{enumerate}
\item[ {\bf (i)}]  $$ \int_0^{\infty} \int_{-\infty}^{vy-\lambda}
\rho'(u+c_0v + g(y) \,v) \;\, v^{n} \, f(u^2+v^2) \; du \;dv\, \leq 0\,, $$
or
\item[ {\bf (ii)}]  $\rho$ is convex, $g \leq g_{\pi_0}$ and
$\psi_{\rho}(\lambda,y) \leq 0$, where
$$\psi_{\rho}(\lambda,y)= \int_0^{\infty}  \int_{-\infty}^{vy-\lambda}
\rho'(u+c_0v + g_{\pi_0}(y) \,v) \;\, v^{n} \, f(u^2+v^2) \; du \;dv\,.$$
\end{enumerate}
\end{theorem}
{\bf  Proof.}  With $\rho'(\cdot)$ increasing by the assumption of convexity, condition {\it \bf (ii)} implies
{\it \bf (i)} so that we only need to establish the sufficiency of {\it \bf (i)}.
Following Kubokawa (1994), write for $\delta_g(X,S) \in C$, 
\begin{align*}
 \,  \rho({ x+c_0s -\mu \over \sigma} ) -
\rho({ x+c_0s + g(\frac{x}{s}) s  -\mu \over \sigma} ) 
& = \rho({ x+c_0s + g(y) s  -\mu \over \sigma} ) |_{y=x/s}^{y = 
\infty} \\
& =  \int_{x/s}^{\infty} \frac{s}{\sigma} \rho'({ x+c_0s + g(y) s -\mu \over
\sigma} )
g'(y) dy \,.
\end{align*}
Now, use the above expression for the difference in losses to write the
difference in risks at $\theta=(\mu,\sigma)$ as:
\begin{align}
\nonumber & \Delta_g(\theta)\\ \nonumber &=  R(\theta, \delta_0) -
R(\theta, \delta_g) \\
\nonumber \; &=  \frac{1}{\sigma} \int_0^{\infty} s \int_{-\infty}^{\infty}
 \{\int_{x/s}^{\infty} g'(y)  \rho'({ x+c_0 s + g(y) s -\mu \over \sigma} ) \,
dy  \} f_{X,S}(x,s) \,dx ds \\
\label{b} \; &= \nonumber \int_{\{g'(y) < 0\}} g'(y)  \{ \int_{0}^{\infty}
\int_{-\infty}^{sy} \rho'({ x+c_0 s + g(y) s -\mu \over \sigma} ) 
 \frac{s^n}{\sigma^{n+2}} f(\frac{(x-\mu)^2+s^2}{\sigma^2})
\, dx ds  \} \,dy,\\&
\end{align}
since $g' \leq 0$ a.e.  
Now, the difference in
risks $\Delta_g(\theta)$ will be nonnegative for all $\theta \in \Theta$ as long as for all $y \in \mathbb{R}$ such that $g'(y)<0$, $\mu \geq 0,
\sigma >0$, the bracketed term in (\ref{b}) is less than or equal than $0$, which is equivalent to {\it \bf (i)}
with the change of variables $(x,s) \rightarrow (u= \frac{x-\mu}{\sigma}, v=\frac{s}{\sigma})$.  \qed

\begin{remark}
\label{expand}
Notice that $\psi_{\rho}(0,y)=B_{n}(y,g_{\pi_0}(y))=0$ for all $y \in \mathbb{R}$ by virtue of the definition of $g_{\pi_0}$ in (\ref{deltapi*def}).
Therefore, the risks of $\delta_{\pi_0}$ and $\delta_0$ match at the boundary of $\Theta$ where $\mu=0$, $\sigma>0$.
Moreover, if $\delta_g $ expands more that $\delta_{\pi_0}$ (whether or not $\delta_g \in C$), then the risk at the boundary of $\delta_g$ will exceed that of 
$\delta_0$, hence giving a condition for non-minimaxity.  This is so given that 
\begin{align*}R((0,\sigma), \delta_g) &=E_{(0,1)}(\rho(\delta_g(X,S)))\\
 &  > E_{(0,1)}(\rho(\delta_{\pi_0}(X,S)))\\&=
R((0,\sigma), \delta_{\pi_0})\\&=R((0,\sigma), \delta_0), \end{align*}
since $\delta_{\pi_0}(X,S) \geq 0$ with probability one, and $\rho$ is increasing on $(0,\infty)$.  
As a consequence of the above, and of Lemma \ref{l} and Corollary \ref{c}, we have the following non-minimaxity result.

\begin{corollary}
\label{non-minimaxity}
For estimating $\mu$ in (\ref{XU}) or (\ref{reduced}) with $\mu \geq 0, \sigma >0$, the generalized Bayes estimators
$\delta_{\pi_l}$ with $-(n-1) \leq l <0$ are not minimax whenever {\bf (a)} $f$ is normal and $\rho$ is even and convex, or
whenever {\bf (b)} $f$ satisfies assumption (\ref{assumptions}) and the loss is invariant $L^p$ with $p > 0$.
\end{corollary}

Analogously, we point out that  $\delta_{\pi_0}$ does not dominate
any other minimax estimator $\delta_g \in C$ taking nonnegative values and satisfying {\it \bf (ii)} of Theorem \ref{main} since such $\delta_g$'s shrink
$\delta_{\pi_0}$ and $R((0,\sigma), \delta_g)=E_{(0,1)}(\rho(\delta_g(X,S))) < E_{(0,1)}(\rho(\delta_{\pi_0}(X,S)))=
R((0,\sigma), \delta_{\pi_0})\,. $
\end{remark}

\begin{remark}
\label{tr}
A plausible alternative to the MRE estimator $\delta_0$ is, of course, its truncation $\delta_{0}^{T}(X,S)=\max(0,\delta_0(X,S))$.
Clearly $\delta_{0}^{T}$ improves upon $\delta_0$ for bowl shaped $\rho$, since for all $\mu \geq 0, \sigma >0$, $\rho(\frac{\delta_{0}^{T}(x,s) - \mu}{\sigma})
\leq \rho(\frac{\delta_{0}(x,s) - \mu}{\sigma})$ for all $(x,s) \in \mathbb{R} \times \mathbb{R}^+$, with strict inequality occurring with positive probability.  Moreover, the estimator $\delta_{0}^{T}$
belongs to the class $C$ with $g_0^T(y)=\max(0,-y-c_0)$, and satisfies condition {\it \bf (i)} of Theorem \ref{main} with 
$\{y: (g_0^{T})'(y) < 0\}=(-\infty,-c_0)$ since 
\begin{align*}
&\int_0^{\infty} \int_{-\infty}^{vy-\lambda}
\rho'(u+c_0v + g_0^{T}(y) \,v) \; v^n \, f(u^2+v^2) \; du \;dv\,\\ 
&= \int_0^{\infty} \int_{-\infty}^{vy-\lambda}
\rho'(u-vy) \; v^n \, f(u^2+v^2) \; du \;dv\ \\
\, 
& \leq  \int_0^{\infty} \int_{-\infty}^{vy-\lambda}
\rho'(-\lambda) \; v^n \, f(u^2+v^2) \; du \;dv \,\\
\, 
& \leq  0\,,
\end{align*}
for all $\lambda \geq 0$.  Finally, the observations of Remark \ref{expand} apply to $\delta_{0}^{T}$, with $\delta_{0}^{T}$ a shrinker of $\delta_{\pi_0}$, and
$\delta_{\pi_0}$ not dominating $\delta_{0}^{T}$.

\end{remark}

With Theorem \ref{main}, our attention focuses on the quantity $\psi_{\rho}(\lambda,y)$ and testing the condition
$\psi_{\rho}(\cdot,\cdot) \leq 0$ on $\mathbb{R}^+ \times \mathbb{R}$ for various choices of $\rho$.  Now, since 
\begin{equation}
\label{0}
\psi_{\rho}(0,y) =0, \; \hbox{and} \lim_{\lambda \to \infty} \psi_{\rho}(\lambda,y)=0 \; \hbox{for all} \; y \in \mathbb{R},
\end{equation}
$\psi_{\rho}(\cdot,y)$ cannot be monotone on $[0,\infty)$ for any $\rho$ and $y \in \mathbb{R}$.
We are thus led to analyzing the behaviour of $\frac{\partial}{\partial \lambda} \psi_{\rho}(\lambda,y)$.

\begin{proposition}
\label{main2}
Let $k(y)=y+c_0+g_{\pi_0}(y)$, $f_{\lambda,y}(t)$ be a Lebesgue density on $(0,\infty)$ proportional to 
$$t^{n}\, f\left(\lambda^2(1+y^2) \, \{(t-\frac{y}{1+y^2})^2 + \frac{1}{(1+y^2)^2}\}\right)\,,$$
and
\begin{equation}
\label{drho} D_{\rho}(\lambda,y) = \int_{0}^{1/k(y)} |\rho'(\lambda(tk(y)-1))| \, f_{\lambda,y}(t) \, dt -
\int_{1/k(y)}^{\infty} |\rho'(\lambda(tk(y)-1))| \, f_{\lambda,y}(t) \, dt \,.
\end{equation}
Suppose further that $D_{\rho}(\cdot,y)$ changes signs once from $-$ to $+$
on $[0,\infty)$ for all $y \in \mathbb{R}$.  Then, for estimating $\mu$ in (\ref{XU}) or (\ref{reduced}) under assumption (\ref{assumptions}) with $\mu \geq 0, \sigma >0$, 
\begin{enumerate}
\item[ {\bf (i)}] the generalized Bayes estimator $\delta_{\pi_0}$ is minimax, 
for strictly bowled shaped loss $\rho(\frac{d-\mu}{\sigma})$ as long as $\delta_{\pi_0} \in C$;

\item[ {\bf (ii)}]  for $\delta_g \in C$, the condition $g \leq g_{\pi_0}$ is sufficient for $\delta_g$ to be minimax
under convex loss $\rho(\frac{d-\mu}{\sigma})$.
\end{enumerate}
\end{proposition}
{\bf Proof.}  We have 
\begin{align}
\nonumber \frac{\partial}{\partial \lambda} \psi_{\rho}(\lambda,y) &= - \int_{0}^{\infty} \, \rho'(vk(y)-\lambda) \, v^n \, f((vy-\lambda)^2 + v^2) \,dv \\
\label{at0} \, & \propto   - \lambda^{n+1} \int_{0}^{\infty} \, \rho'(\lambda t k(y)-\lambda) \, f_{\lambda,y}(t) \, dt \\
\label{psiD} \, & \propto  D_{\rho}(\lambda,y)\,.
\end{align}
Therefore, under the given assumptions on the sign changes of $D_{\rho}(\cdot,y)$, we infer that, for all
 $y \in \mathbb{R}$,  
$\psi_{\rho}(\lambda,y)$ decreases, then increases as $\lambda$ varies on $[0,\infty)$.  Finally, the result follows from Theorem \ref{main}
and property (\ref{0}).  \qed

\begin{remark}
\label{remarks}
\begin{enumerate}
\item[ {\bf (i)}]
From (\ref{at0}), note that $$\frac{\partial}{\partial \lambda} \psi_{\rho}(\lambda,y)|_{\lambda=0^+}= - \int_{0}^{\infty} \, \rho'(vk(y)) \, v^n \, f(v^2(y^2+1)) \,dv \leq 0,$$ since $k(\cdot) \geq 0$ from (\ref{positive}), and $\rho'(\cdot) \geq 0$ on $[0,\infty)$.  Hence, Proposition \ref{main2}'s sign change assumption on $D_{\rho}(\cdot,y)$ is consistent, for any strictly bowled-shaped $\rho$, with the behaviour of $\psi_{\rho}(\lambda,y)$ for $\lambda$ near $0$.  
\item[ {\bf (ii)}]
Turning to the families of densities $\{f_{\lambda,y}(\cdot), \lambda \in [0,\infty), y \in \mathbb{R}\}$, they can be shown for $y \leq 0$ to possess a decreasing 
monotone likelihood ratio (mlr) in $T$, or equivalently in $W=(T- \frac{y}{y^2+1})^2$, with $\lambda$ viewed as the parameter.  Indeed, for $\lambda_1 > \lambda_0 \geq 0$,
setting $\alpha_i=\lambda_i^2 (y^2+1)$ and $\epsilon=(y^2+1)^{-2}$, we have 
$$ \frac{f_{\lambda_1,y}(t)}{f_{\lambda_0,y}(t)} \, \propto \, \frac{f(\alpha_1 (w + \epsilon))}{f(\alpha_0 (w + \epsilon))}$$
which decreases in $w$, $w>\frac{y^2}{(y^2+1)^2}$, given assumption (\ref{assumptions}).

\item[ {\bf (iii)}]
In the normal case, the densities $f_{\lambda,y}(\cdot)$ may 
described as weighted (by the factor $t^{n}$) positively truncated $N(y/(1+y^2),1/(\lambda^2(1+y^2)))$ densities.  They have been recently studied in related work of Marchand, Jafari Jozani, and Tripathi (2011) where quantiles are estimated under the restriction $\mu \geq 0$. 
\end{enumerate}
\end{remark}

We conclude this section with a very useful technical result.

\begin{lemma}
\label{k}  Let $k(y)=y+c_0+g_{\pi_0}(y)$ as in Proposition \ref{main2} and let $\rho$ be either an even function, or more generally satisfy $|\rho'(-u)| \leq |\rho'(u)|$ for all $u>0$.  Then we have $\frac{1}{k(y)} > \max\{0,\frac{y}{1+y^2}\}$.
\end{lemma}
{\bf Proof.}  The positivity of $k(y), y \in \mathbb{R}$ follows from (\ref{positive}).  To establish that
$\frac{1}{k(y)} > \frac{y}{1+y^2}$, we assume the contrary and show that this would imply
$D_{\rho}(\lambda,y) \leq 0$ for all $\lambda \geq 0$ which is not possible given (\ref{at0}) and (\ref{psiD}).
Indeed, we would have, for all $\lambda \geq 0, y \in \mathbb{R}$, under the given assumption on $\rho$
\begin{align*}
&D_{\rho}(\lambda,y) \\
&\leq \int_{0}^{1/k(y)} |\rho'(\lambda k(y)(t-\frac{1}{k(y)}))| \, f_{\lambda,y}(t) \, dt -
\int_{1/k(y)}^{2/k(y)}  |\rho'(\lambda k(y)(t-\frac{1}{k(y)}))| \, f_{\lambda,y}(t) \, dt \\
\, 
& \leq   \int_{0}^{1/k(y)} |\rho'(\lambda k(y)(t-\frac{1}{k(y)}))| \;\, (f_{\lambda,y}(t) - 
f_{\lambda,y}(\frac{2}{k(y)}-t) )\, dt \\
& \leq  0\,,
\end{align*}
given that $f_{\lambda,y}(t) \leq f_{\lambda,y}(\frac{2}{k(y)}-t)$ for all $t \in (0, 1/k(y))$ whenever $\frac{1}{k(y)} \leq \frac{y}{1+y^2}$.
\qed 

The inequality $\frac{1}{k(y)} > \max\{0,\frac{y}{1+y^2}\}$ will be exploited as a technical result, but it also provides an interesting upper bound for the generalized Bayes estimator $\delta_{\pi_0}$, namely 
$$ \delta_{\pi_0}(x,s) = s \, k(\frac{x}{s}) < x + \frac{s^2}{x}, \; \hbox{for} \; x>0,$$ 
applicable to all pairs $(f,\rho)$ for which $\delta_{\pi_0}$ exists, with $f$ satisfying (\ref{assumptions}), $\rho$ satisfying the conditions of Lemma \ref{k}.

\section{Minimax results for the normal case}

Here is a minimax result applicable in the normal case, to Bayes estimators $\delta_{\pi_l}$, and for general convex losses that are either
even functions or, more generally, that penalize the rate of over-estimation more sharply than the rate of underestimation in the sense 
\begin{equation}
\label{asymmetrycondition}
|\rho'(-u)| \leq \rho'(u)\,, \; \; \hbox{for all} \; \; u \geq 0 \,.  
\end{equation}

\begin{theorem}
\label{normalcase}
For estimating $\mu$ in the normal case in (\ref{model}) with $\mu \geq 0, \sigma >0$ under convex $\rho$ in (\ref{loss}),
\begin{enumerate}
\item[ {\bf (a)}] the condition $g \leq g_{\pi_0}$ suffices for an estimator $\delta_g \in C$ to 
be minimax in cases where $\rho$ satisfies condition (\ref{asymmetrycondition});

\item[ {\bf (b)}] such minimax estimators include the generalized Bayes estimator $\delta_{\pi_0}$ under losses $\rho$ satisfying (\ref{asymmetrycondition}),
and all $\delta_{\pi_l}$ with $l >0$ when $\rho$ is even.
\end{enumerate}
\end{theorem}
{\bf Proof.}   Given that $\delta_{\pi_l} \in C$ for $l=0$, and for $l >0$ when $\rho$ is even by virtue of Lemma \ref{gpi*}, 
the first part of (b) is simply a restatement of (a) for the Bayes estimator $\delta_{\pi_0}$, while the part relating to
$\delta_{\pi_l}$ with $l >0$ follows also from (a) and Lemma \ref{l}.  The rest of the proof concerns part (a) and we apply 
Proposition \ref{main2}.  From (\ref{drho}), we have with the change of variables $u=\lambda(tk(y)-1)$:
$$  D_{\rho}(\lambda, y) = \frac{1}{\lambda k(y)} E[-\rho'(U)]\,,$$
where $U$ has density proportional to 
\begin{equation}
\label{densityU}
f_{\lambda, y} \left((\frac{u}{\lambda}+1) \frac{1}{k(y)}  \right) 1_{(-\lambda, \infty)}(u)\,.
\end{equation}
Since $-\rho'$ changes signs once from $+$ to $-$ on $\mathbb{R}$, a decreasing in $u$ monotone likelihood ratio
property of the densities in (\ref{densityU}) with respect to the parameter $\lambda$ 
will suffice to establish that $D_{\rho}(\lambda, y)$ changes signs from $-$ to $+$ on $[0,\infty)$ as a function of $\lambda \geq 0$
and permit us to apply Proposition \ref{main2}\footnote{It is interesting to point out that the arguments here apply as well to strictly-bowled shaped losses.  As well, only a stochastic increasing property for the densities $f$ is required. The monotonicity of $g_{\pi_0}$ however is guaranteed by the convexity of $\rho$ (Lemma \ref{gpi*}), which is assumed here.}.   Now, the densities in (\ref{densityU}) may be written as 
$$ h_{\lambda}(u) \propto (\frac{u+\lambda}{\lambda})^n \, f\left(c(u+\lambda b)^2 + d\lambda^2  \right) 1_{(-\lambda, \infty)}(u)\,,$$
with $c=(1+y^2)/k^2(y)$, $b=1- \left(yk(y)/(1+y^2)\right)$, and $d=(1+y^2)^{-1}$.  Notice that we have $c>0$ by virtue of (\ref{positive}), and $b>0$ by 
assumption (\ref{asymmetrycondition}) and Lemma \ref{k}.  Finally, in the normal case with $f(t) = (2\pi)^{(n+1)/2}e^{-t/2}$, tha ratio 
$\frac{h_{\lambda_1}(u)}{h_{\lambda_0}(u)}$ is, for $\lambda_1 > \lambda_0 \geq 0$, undetermined for $u \leq -\lambda_1$, equal to $+\infty$ for 
$u \in (-\lambda_1, -\lambda_0]$, and otherwise proportional to 
$$ \left(\frac{u+\lambda_1}{u+\lambda_0}\right)^n \, e^{-bcu(\lambda_1 - \lambda_0)},$$
which is indeed decreasing in $u$ for $u > -\lambda_0$, and which establishes the result. \qed

The normal case minimax results of Theorem \ref{normalcase} in part {\bf (a)}, and applicable to the generalized Bayes estimator $\delta_{\pi_0}$, were previously obtained
for the specific case of scale invariant $L^2$ loss by Kubokawa (2004).  He works directly with the Bayes estimator in Example 1 to derive the key required analytical properties, namely the monotonicity of $g_{\pi_0}$ in Lemma \ref{gpi*} and inequality {\bf (i)} of Theorem \ref{main}.  With Kubokawa's analysis specific to scale invariant $L^2$ loss, the normal model and the estimator $\delta_{\pi_0}$, our unified development above contrasts and provides extensions with respect to the loss and the prior.  In the next section, we give extensions with respect to the model for scale invariant $L^p$ losses and asymmetric versions.

\section{Minimax results for scale invariant $L^p$ losses and their asymmetric versions}

The minimax results of this section are applicable for the wider class of models, or choices of $f$, in (\ref{XU}) with assumptions (\ref{assumptions}).
As well, these findings concern scale invariant $L^{p}$ losses $|\frac{d-\mu}{\sigma}|^p$, $p>0$, and the more 
general $\rho_{c_1,c_2}$ in (\ref{lp}) with $c_2 \geq c_1$.   For these losses, 
(\ref{drho}) reduces to $D_{\rho}(\lambda,y)=p\lambda^{p-1} E_{\lambda}[\,g_y(T)\,]$, with $T \sim f_{\lambda,y}$ and $g_y(t)=c_1 \, (1-tk(y))^{p-1} 1_{(0,1/k(y)}(t)
-c_2 \, (tk(y)-1)^{p-1} 1_{(1/k(y),\infty)}(t)$.  With this representation, observe that $g_y(\cdot)$ changes sign once on $(0,\infty)$ from $+$ to $-$,
so that $E_{\lambda}[\,g_y(T)\,]$ changes signs from $-$ to $+$ as $\lambda$ varies on $[0,\infty)$, in view of sign change properties and the mlr property of Remark \ref{remarks} {\bf (ii)}.  Therefore, $D_{\rho}(\lambda,y)$ varies indeed, as a function of $\lambda \in [0,\infty)$ from $-$ to $+$ as prescribed in Proposition \ref{main2} for $y \leq 0$ and losses $\rho_{c_1,c_2}$.  For $y>0$ however, the situation is more delicate.  We continue with the non-convex case with $p \in (0,1)$, and this will be followed by the convex case with $p \geq 1$.

\begin{theorem}
\label{0<p<1}
For estimating $\mu$ in (\ref{XU}) or (\ref{reduced}) with $\mu \geq 0, \sigma >0$ under scale invariant $L^p$ loss $|\frac{d-\mu}{\sigma}|^p$ with
$p \in (0,1)$, the generalized Bayes estimator $\delta_{\pi_0}$ is minimax.
\end{theorem}
{\bf Proof.} With $\delta_{\pi_0} \in C$ by virtue of Lemma \ref{gpi*}, we seek to apply part {\bf (i)} of Proposition \ref{main2} to show that
$\delta_{\pi_0}$ is minimax.  For $\rho(t)=|t|^p$ with $p>0$, we reexpress (\ref{drho}) as
\begin{equation}
\nonumber
D_{\rho}(\lambda, y) \propto E[A(T)\, B(T)] = E[G(S)]\,, 
\end{equation}
with 
\begin{align*}
G(s)=E[A(T)\, B(T)|S=s],\, A(t)= t^n \, |t-\frac{1}{k(y)}|^{p-1},\, B(t)= -1 + 2 I_{(0,\frac{1}{k(y)}]}(t), \end{align*}
\begin{align*}
T \sim f\left(\lambda^2(1+y^2) \, \{(t-\frac{y}{1+y^2})^2 + \frac{1}{(1+y^2)^2}\}\right), 
\end{align*}
and 
\begin{align*}
 S \overset{d}{=} (T-\frac{y}{1+y^2})^2 + \frac{1}{(1+y^2)^2}\,.
\end{align*}
Given assumption (\ref{assumptions}), the family of densities of $S$ are seen to have a decreasing monotone likelihood ratio in $S$ with parameter $\lambda^2(1+y^2)$.
So, in accordance with Karlin's sign change analysis, to prove the result, it will suffice to show that
\begin{equation}
\label{signchangeG}
G(s) \;\,\hbox{changes signs once from} \;\, + \; \hbox{to} \; - \; \hbox{as} \; s \; \hbox{varies on} \; (\frac{1}{(1+y^2)^2},\infty)
\end{equation}
to establish that $D_{\rho}(\lambda, y)$ changes signs as prescribed by Proposition \ref{main2}.  We proceed by treating separately the cases: 
{\bf (i)} $\frac{1}{k(y)} - \frac{y}{1+y^2} \geq \frac{y}{1+y^2}$  and {\bf (ii)} $0 \leq \frac{1}{k(y)} - \frac{y}{1+y^2} \leq \frac{y}{1+y^2}$.
Here, we have made use of Lemma \ref{k} to discount the remaining possibility $\frac{1}{k(y)} - \frac{y}{1+y^2} < 0$.
\begin{enumerate}
\item[ Case {\bf (i)}:]  Set $s_0= (\frac{1}{k(y)} - \frac{y}{1+y^2})^2 + \frac{1}{(1+y^2)^2}$.  Observe that, whenever $s \geq s_0$,
$P(T \geq \frac{1}{k(y)}|S=s)=1$ implying $P(B(T)=-1|S=s)=1$ and $G(s) \leq 0$.  Similarly, if $s < s_0$, then $P(B(T)=1|S=s)=1$ and $G(s) \geq 0$.
Hence, the above establishes (\ref{signchangeG}) for case {\bf (i)}.

\item[ Case {\bf (ii)}:]  Here, we set $s_1= \frac{y^2}{(1+y^2)^2} + \frac{1}{1+y^2}$, so that $s_0 \leq s_1$.  As in {\bf(i)}, we verify that $G(s) \geq 0$
for $s \leq s_0$, and $G(s) \leq 0$ for $s \geq s_1$.  Finally, for $s \in (s_0, s_1)$, the conditional distribution of $T|S=s$ is a two-point uniform discrete
distribution on $\{t_1, t_2 \}$, with $t_1 = \frac{y}{1+y^2} + \Delta$,  $ t_2=  \frac{y}{1+y^2} - \Delta$ and $\Delta= \sqrt{s-\frac{1}{(1+y^2)^2}}\,$.
\end{enumerate}
We hence obtain
\begin{align*}
G(s) &= \frac{1}{2} (A(t_2) - A(t_1)) \\
\, & =  \frac{1}{2} \left[ t_2^n \; |t_2 - \frac{1}{k(y)}|^{p-1} - t_1^n \; |t_1 - \frac{1}{k(y)}|^{p-1}    \right]  < 0\,, 
\end{align*}
since $ t_2 < t_1$, $n >1$; $|t_1 - \frac{1}{k(y)}|=  \frac{y}{1+y^2} - \frac{1}{k(y)} + \Delta < - \frac{y}{1+y^2} + \frac{1}{k(y)} + \Delta = |t_2 - \frac{1}{k(y)}|$, and $p-1 <0$. Hence, the above establishes (\ref{signchangeG}) for case {\bf (ii)} and completes the proof. \qed

\begin{theorem}
\label{main3}
For estimating $\mu$ in (\ref{XU}) or (\ref{reduced}) under assumptions (\ref{assumptions}), with $\mu \geq 0, \sigma >0$ and with loss
$\rho_{c_1, c_2}$, $p \geq 1$ and $c_2 \geq c_1$,
\begin{enumerate}
\item[ {\bf (a)}] the condition $g \leq g_{\pi_0}$ suffices for an estimator $\delta_g \in C$ to be minimax;

\item[ {\bf (b)}]  such minimax estimators include the generalized Bayes estimator $\delta_{\pi_0}$, as well as all generalized
Bayes estimators $\delta_{\pi_l}$, $l>0$ for the symmetric case $c_1=c_2$.
\end{enumerate}
\end{theorem}

{\bf Proof.} For losses $\rho$ as in (\ref{lp}), we may write 
\begin{equation}
\nonumber
D_{\rho}(\lambda,y)= p \lambda^{p-1} \lbrace c_1 \int_0 ^{\frac{1}{k(y)}} h_{\lambda,y}(w) \, dw - 
c_2 \int_{\frac{1}{k(y)}}^{\infty} h_{\lambda,y}(w) \, dw \, \rbrace,
\end{equation}
with $h_{\lambda,y}(\cdot)$ a probability density function on $(0,\infty)$ proportional to $ |wk(y)-1|^{p-1} \,f_{\lambda,y}(w)$.
From this, we see that $D_{\rho}(\lambda,y)$ is positive iff $P_{\lambda}(W>1/k(y)) < c_1/(c_1+c_2)$, 
where $W$ is a random variable with pdf $h_{\lambda,y}$.  We show below in Section \ref{q} of the Appendix that, whenever $\frac{1}{k(y)} > \frac{y}{1+y^2}$, the quantity
$P_{\lambda}(W>1/k(y))$ decreases in $\lambda$ on $[0,\infty)$, which means that $D_{\rho}(\cdot,y)$ changes signs from $-$ to $+$ on $[0,\infty)$.
The result then follows from Proposition \ref{main2} and Lemma \ref{k}.  \qed

\section{Concluding Remarks }

We have considered the problem of estimating a lower bounded location parameter for a wide array of spherically symmetric location-scale models with a residual vector as represented in model (\ref{XU}), with unknown scale, and under scale invariant loss as given by (\ref{loss}).  With a relative paucity of findings for such problems when the scale parameter is unknown, we have established the minimaxity of the generalized Bayes estimator $\delta_{\pi_0}$ for normal models and convex loss, as well for more general models and scale invariant $L^p$ loss and asymmetric versions $\rho_{c_1,c_2}$ given in
(\ref{lp}).  Moreover, we have shown the role of $\delta_{\pi_0}$ to be pivotal, in the sense that it provides an upper threshold condition necessary for the minimaxity 
of many estimators.   Other minimax estimators are also obtained, including generalized Bayes estimators $\delta_{\pi_l}$ when $l>0$ and the loss is convex and even in the above situations.   The results represent extensions of Kubokawa's results (2004) applicable to scale invariant $L^2$ loss.  Much of the treatment is unified and exploits general features of the model and the loss with incisive analysis and novel representations.    Various other observations are given, including the robustness
of the Bayes estimator $\delta_{\pi_l}$ with respect to the choice of $f$ in model (\ref{XU}).

As illustrated by Marchand, Jafari Jozani and Tripathi (2011), the normal case improvements provided for scale invariant $L^2$ loss yield applications for two-sample problems where $Y_i \sim N(\mu_i, \sigma^2); i=1,2$ with unknown $\mu_1, \mu_2, \sigma^2,$
where the objective is to estimate $\mu_1$ (or $\mu_2$) with the additional information of the ordering $\mu_1 \leq \mu_2$.
Despite the advances presented here, minimax extensions to other strictly bowled-shaped losses, although plausible, are still lacking.  Furthermore, numerous questions remain unanswered such as the admissibility of the above minimax estimators,
the investigation of wider classes of Bayes estimators for minimaxity, and related tests of minimaxity for multivariate location-scale problems with order restrictions.

\section{Appendix}

\subsection{\bf Proof of Lemma \ref{mixtures}}
We have 
\begin{align}
\label{t} 
\nonumber \frac{t\, f'(t)}{f(t)}  &=  \frac{\int_{0}^{\infty} t v^2 \, f_0'(tv) \, h(v) \, dv }
{\int_{0}^{\infty} v \, f_0(tv) \, h(v) \, dv} 
\, \\ \nonumber &=  \frac{\int_{0}^{\infty} z \, f_0'(z) \, h(z/\,t) \, dz}
{\int_{0}^{\infty}  f_0(z) \, h(z/\,t) \, dz} \\&=  E_t[\frac{Z\, f_0'(Z)}{f_0(Z)}]\,,
\end{align}
where $Z$ has density proportional to $f_0(z) \, h(z/\,t)$ on $\mathbb{R}^+$.
Now, observe that this scale family of densities for $Z$ has increasing monotone likelihood ratio in $Z$, with parameter $t$,
as a consequence of assumption (\ref{assumptions}) for $h$.  Finally, the result follows from representation (\ref{t}) with this monotone likelihood ratio 
and since $\frac{zf_0'(z)}{f_0(z)}$ decreases in $z$ by assumption (\ref{assumptions}) for $f_0$.  \qed

\subsection{\bf Proof of Lemma \ref{l}}

We fix $y \in \mathbb{R}$ throughout and set $c_0=0$ given that $\rho$ is assumed even.  First, observe that by differentiating (\ref{deltapi*def}) for the normal case with $f(u^2+v^2) = \phi(u) h(v)$ in (\ref{densities}), we have $\frac{\partial}{\partial \;l} B_{n+l}(y, g_{\pi_l}(y) = 0$ which implies
\small
\begin{align*}
& \int_0^{\infty} [\frac{\partial}{\partial \; l} \left\{\int_{-\infty}^{vy} \rho'(u+ g_{\pi_l}(y) v) \phi(u) du \right\}  \\
&+ \left\{ \int_{-\infty}^{vy} \rho'(u+ g_{\pi_l}(y) v) \phi(u) du \right\} \log(v)] v^{n+l} h(v) dv =0.\end{align*}
\normalsize
Hence, given that $\rho'$ is increasing, to show that $g_{\pi_l}(y)$ decreases in $l$, it will suffice to show that $I \geq 0$, where
\begin{equation}
I = \int_0^{\infty} (\log v) A_y(v) v^{n+l} h(v) dv, \;\;
\hbox{and} \;\; A_y(v)= \int_{-\infty}^{vy} \rho'(u+ g_{\pi_l}(y) v) \phi(u) du\,. 
\end{equation}
Now, we will show below that  
\begin{equation}
\label{ayv}
A_y(v) \;\; \hbox{changes signs once as a function of } \;\; v \;\; \hbox{from} - \;\;\hbox{to} \;\; +\,.
\end{equation}
Applying Lemma \ref{signchange}, which is stated in the Appendix, with $\xi \sim \xi^{n+l} \,h(\xi) \,1_{(0,\infty)}(\xi)$, $r(\xi)=\log(\xi)$, and
$s(\xi)=A_y(\xi)$, we infer that $I \geq 0$ since $E[A_y(\xi)]=0$ given the definition of $g_{\pi_l}(y)$ in (\ref{deltapi*def}).
There remains to establish (\ref{ayv}), which we proceed to do separating the cases: {\bf (i)} $y \leq 0$ and {\bf (ii)} $y>0$. 
\begin{enumerate}
\item[ {\bf (i)}] Case $y<0$.  Let $v_0$ be such that $A_y(v_0)=0$.  Such a value exists since the average value of $A_y(\xi)$ under the above density for $\xi$ is equal to $0$.  For $\epsilon \geq 0$, we have
\begin{align*}
&\frac{A_y(v_0+\epsilon)}{\Phi((v_0 + \epsilon)y)}\\ & =  \int_0^{(v_0 + \epsilon)y} \rho'(u-\epsilon y +
\epsilon(y+ g_{\pi_l}(y)) +  g_{\pi_l}(y) v_0) \frac{\phi(u)}{\Phi((v_0 + \epsilon)y)} \, du \\
\, & \geq  \int_0^{v_0 y} \rho'(u' +  g_{\pi_l}(y) v_0) \frac{\phi(u' + \epsilon y)}{\Phi((v_0 + \epsilon)y)} \, du' \\
\, & =  C_y(v_0, \epsilon) \; \; (\hbox{say})\,,
\end{align*}
with equality if and only if $\epsilon=0$, given (\ref{positive}) and since $\rho'$ is increasing.  Now, observe that the ratio
of densities $\frac{(\phi(u'+\epsilon y)/\Phi((v_0 + \epsilon)y))}{\phi(u')/\Phi(v_0 + y)}$ is increasing in $u'$, for $u' \in (-\infty, v_0y)$
and $\epsilon y < 0$.  Hence, this monotone likelihood ratio property implies that $\frac{A_y(v_0+\epsilon)}{\Phi((v_0 + \epsilon)y)} 
\geq C_y(v_0,\epsilon) \geq C_y(v_0,0)=0$, for $\epsilon >0$ and $y<0$,
yielding (\ref{ayv}) for $y<0$.

\item[ {\bf (ii)}] Case $y \geq 0$.  As in {\bf (i)}, let $v_0$ be such that $A_y(v_0)=0$.  Using this, as well as property (\ref{positive}), the nonnegativity of $g_{\pi_l}(\cdot)$ (Lemma \ref{gpi*}) (since $c_0=0$), and the convexity of $\rho$, we have for $\epsilon >0$:
\begin{align*}
A_y(v_0+\epsilon) & =  \int_{-\infty}^{(v_0+\epsilon) y} \rho'\{u+(g_{\pi_l}(y)) (v_0+\epsilon) \} \, \phi(u) du\, \\
\, & \geq  A_y(v_0) + \int_{v_0 y}^{(v_0+\epsilon) y} \rho'\{u+(g_{\pi_l}(y)) v_0 \} \, \phi(u) du\, \\
& \geq 0.  \;\;\;\;\;\;\;\;\;\;\;\;\;\;\;\;\;\;\;\;\;\;\;\;   \qed 
\end{align*}  
\end{enumerate}

\subsection{Lemma used within the proof of Lemma \ref{l}.}

The following result is well known and its proof is left to the reader.

\begin{lemma}
\label{signchange}
Let $\xi$ be a continuous random variable, let $r(\cdot)$ be a continuous and increasing function on the support
of $\xi$, and let $s(\cdot)$ be a continuous function which changes signs once from $-$ to $+$ at $s_0$ on the support of
$\xi$.  Then, we have $E[\,r(\xi) \, s(\xi)\,]\,\geq \, r(s_0) \, E[s(\xi)]$, and, in particular if $E[s(\xi)]=0$, then 
$E[\,r(\xi) \, s(\xi)\,]\, \geq 0$.
\end{lemma}

\subsection{\bf Theorem \ref{main3}: proof of a monotonocity property for $P_{\lambda}(W>1/k(y))$} \label{q}  We wish to show that  
\begin{equation}
\label{d}
P_{\lambda}(W>\frac{1}{k(y)}) \; \, \hbox{decreases in } \,\alpha \; \hbox{whenever} \,\; \frac{1}{k(y)} > a\,,
\end{equation}
with $W$ having pdf $h_{\lambda,y}(w)$ on $(0,\infty)$ proportional to $ |wk(y)-1|^{p-1} \,w^{n} \,f(\alpha \{(w-a)^2 + \epsilon\}$,
$a=\frac{y}{1+y^2}$, $\alpha=\lambda^2(1+y^2)$, and $\epsilon=(y^2+1)^{-2}$.  We have 
\begin{align} & \nonumber \frac{\partial}{\partial \lambda} P_{\lambda}(W>\frac{1}{k(y)}) \\&=\nonumber  2\lambda (1+y^2)
\frac{\partial}{\partial \alpha} {\int_{\frac{1}{k(y)}}^{\infty} |wk(y)-1|^{p-1} \,w^{n} \, f(\alpha \{(w-a)^2 + \epsilon\} \, dw \over 
 \int_{0}^{\infty} |wk(y)-1|^{p-1} \, w^{n} \, f(\alpha \{(w-a)^2 + \epsilon\} \, dw} \leq 0 \\
\label{keyappendix}  \, & \Longleftrightarrow  E[\gamma(\alpha\{(W-a)^2 + \epsilon \}) |W> \frac{1}{k(y)}] \geq
E[\gamma(\alpha\{(W-a)^2 + \epsilon \})],  
\end{align} 
under pdf $h_{\lambda,y}$, with $\gamma(t)=t \frac{|f'(t)|}{f(t)}$.
Taken together, the following hence form a sufficient condition for (\ref{keyappendix}) to hold:  $$ 
{\bf (i)} E[\gamma(\alpha\{(W-a)^2 + \epsilon \})|W> \frac{1}{k(y)}] \geq E[\gamma(\alpha\{(W-a)^2 + \epsilon \})|W<a]\,, \text{and}
$$ $$ {\bf (ii)} \;
E[\gamma(\alpha\{(W-a)^2 + \epsilon \})|W> \frac{1}{k(y)}] \geq E[\gamma(\alpha\{(W-a)^2 + \epsilon \})|a \leq W< \frac{1}{k(y)}]\,.$$
Condition {\bf (ii)} is immediate since $\gamma(\alpha\{(W-a)^2 + \epsilon \})$ increases in $(W-a)^2$ on $(a,\infty)$ by assumption (\ref{assumptions}).  For {\bf (i)}, set $Z=|W-a|$ so that 
$Z|W>1/k(y)$ has pdf proportional to $|(a+z)k(y)-1|^{p-1} \, (a+z)^{n} \,f(\alpha(z^2+\epsilon)) \, 1_{(1/k(y)-a, \infty)}(z)$, while
$Z|W<a$ has pdf proportional to $|1-(a-z)k(y)|^{p-1} \, (a-z)^{n} \,f(\alpha(z^2+\epsilon)) \, 1_{(0,a)}(z)$.  We thus have the ratio
\begin{eqnarray*}
\frac{f_{Z|W<a}(z) }{f_{Z|W>\frac{1}{k(y)}}(z)}\propto\left\{
                                                              \begin{array}{ll}
                                                                \infty &  \hbox{if $z<a, z<\frac{1}{k(y)}-a$;} \\
                                                                0 & \hbox{if $z\geq a, z>\frac{1}{k(y)}-a$;} \\
                                                                \left(\frac{a-z}{a+z}\right)^{n}\left(\frac{zk(y)+(1-ak(y))}{zk(y)-(1-ak(y))}\right)^{p-1} & \hbox{if $0<z<a, z>\frac{1}{k(y)}-a$.}
                                                              \end{array}
                                                            \right.
\end{eqnarray*}

Since both $\frac{a-z}{a+z}$ and
$\frac{zk(y)+b}{zk(y)-b}$ decrease in $z$ for $z<a$ and $z>1/k(y)-a$, 
with $b=1-ak(y)>0$ (Lemma \ref{k}), we have a decreasing monotone likelihood ratio.  Finally, with $\gamma(\alpha(z^2+\epsilon))$ increasing in $z>0$ by (\ref{assumptions}), condition {\bf (i)} follows and our proof of (\ref{d}) is complete.    \qed

\section*{Acknowledgments}
This work was partially supported by a grant from the Simons Foundation (\#209035) to William Strawderman.
Mohammad Jafari Jozani and \'Eric Marchand gratefully acknowledge the research support
of the Natural Sciences and Engineering Research Council of Canada.



\end{document}